\documentclass{amsart}
 \newtheorem{thm}{Theorem}[section]
 
 \newtheorem{lem}[thm]{Lemma}
 
 \theoremstyle{definition}
 
 \theoremstyle{remark}

 \numberwithin{equation}{section}

\newcommand{\norm}[1]{\left\Vert #1 \right\Vert}
\newcommand{\prob}[1]{ \mathbb{P} \left( #1 \right) }
\newcommand{\dprob}[0]{\ \mathrm{d} \mathbb{P}}
\newcommand{\erw}[1]{ \mathbb{E} \left[ #1 \right] }
\newcommand{\D}[0]{\ \mathrm{d}}

\begin{document}

\title[Khintchine inequality for $S^{N-1}$-variables in $L_{\psi_2}$]
 {The best constant in the Khintchine inequality of the Orlicz space $L_{\psi_2}$ for equidistributed random variables on spheres}

\author[H. Dirksen]{Hauke Dirksen}
\address{%
Department of Mathematics\\
Kiel University}
\email{dirksen@math.uni-kiel.de}
\date{March 8, 2016}

\begin{abstract}
We compute the best constant in the Khintchine inequality for equidistributed random variables on the $N$-sphere in the Orlicz space $L_{\psi_2}$. 
\end{abstract}

\maketitle
\section{Introduction}
The classical Khintchine inequality compares the $L_p$-norm of a sum of Rade\-macher variables with the $\ell_2$-norm of the coefficients of the sum. The computation of the best possible constants has attracted a lot of interest. For the classical case, Haagerup found the best constants  for general $p\in (1,\infty)$ in \cite{Haagerup1982}.
Also Khintchine inequalities for different kinds of random variables were investigated, for example rotationally invariant random vectors in \cite{Koenig2001}.
A second variation of the problem changes the underlying space. The Khintchine inequality in Orlicz spaces has been considered in various cases, the first example is a paper by Rodin and Semyonov \cite{Rodin1975}.

Let $q>0$ and $\psi_q(x):= \exp(x^q)-1$ for $x\in \mathbb{R}$. By $\norm{\cdot}_{\psi_q}$ we denote the norm of the Orlicz space $L_{\psi_q}(\Omega, \Sigma,\mu)$. This is given by
\[
\norm{X}_{\psi_q}:=\inf\{c>0 \mid \erw{\psi_q\left(\frac {\norm{X}} c\right)}\leq 1\},
\]
for $X\in L_{\psi_q}$. By $\norm{\cdot}$ we denote the Euclidean norm.
For $q\leq 2$ one can still compare the $L_{\psi_q}$-norm and the $\ell_2$-norm, see \cite{Ledoux1991}. For $q>2$, Pisier proved that the Lorentz sequence spaces $\ell_{q',\infty}$  $(1/q + 1/q' =1)$, instead of $\ell_2$ come into play, see \cite{Pisier1981}. This fact was already mentioned by Rodin and Semyonov \cite{Rodin1975}.

Here we compute the best constant for the Orlicz space $L_{\psi_2}$ and equidistributed variables on $N$-dimensional spheres. 
We apply the technique from \cite{Peskir1993}. Peskir reduces the case of the Orlicz space to the classical Khintchine inequality in $L_q$. The optimality of the constants from $L_q$ carries over to $L_{\psi_2}$.
The same reduction technique can be used for variables on spheres. K\"onig and Kwapien computed the optimal constants in \cite{Koenig2001}. Again the optimality carries over. In this paper we prove the following result.
\begin{thm}\label{thm1}
Let $X_j$, $j=1,\dots, n$ be an i.i.d. sequence of  equidistributed random variables on the $N$-sphere $S^{N-1}$.
For all $a=a_1,\dots,a_n \in \mathbb{R}$ we have
\[
\norm{\sum_{j=1}^{n} a_j X_j}_{\psi_2}\leq b(N) \bigg(\sum_{j=1}^n a_j^2\bigg)^{\frac 1 2 },
\]
where the constant $b(N):=\sqrt{\frac 2 N}{\sqrt{\frac{1}{1-(\frac 1 2 )^{\frac 2 N }}}}$ is optimal.
\end{thm}
Note that $b(N)$ decreases to $\frac{1}{\sqrt{\ln 2}}$ for $N\to \infty$.
In Section \ref{sec:inequality} we prove that the inequality is true. Therefore we consider the series expansion of the exponential function. Then we apply the Khintchine inequality from \cite{Koenig2001}.
In Section \ref{sec:optimality} we show that the constant $b(N)$ can not be smaller. We show that with $Y_n:=\sum_{j=1}^{n}\frac 1{\sqrt{n}}X_j$ we get asymptotic equality in Theorem \ref{thm1} for $n\to\infty$.
\section{Proof of the inequality}\label{sec:inequality}
Let $C>0$. Applying Beppo-Levi we may interchange the limit and the expected value.
\begin{align}\label{eq:seriesexp}
		&\erw{ \exp \left( \frac{\norm{\sum_{j=1}^n a_j X_j}^2}{C^2 \sum_{j=1}^n \norm{ a_j } ^2} \right)} \nonumber \\
		=&
		\erw{\sum_{k=0}^{\infty} \frac{1}{k!}\frac{1}{C^{2k} \left(\sum_{j=1}^n \norm{ a_j } ^{2}\right)^k } \norm{\sum_{j=1}^n a_j X_j}^{2k} } \nonumber\\		
		=& \sum_{k=0}^{\infty} \frac{1}{k!} \frac{1}{ C^{2k} \left(\sum_{j=1}^n\norm{ a_j } ^{2}\right)^k} \erw{\norm{\sum_{j=1}^n a_j X_j}^{2k} }
\end{align}
Now we apply K\"onig's and Kwapie\'n's Khintchine inequality for variables on the sphere and use the constants for $p=2k$, which gives $\left(\widetilde {b}(2k)\right)^{2k}=\left(\frac 2 N\right) ^k\left( \frac{\Gamma(k+\frac N 2)}{\Gamma \frac N 2}\right)$, see \cite[Theorem 3]{Koenig2001}. We obtain
\[
	\erw{\norm{\sum_{j=1}^n a_j X_j}^{2k}} 
			\leq 
	\left( \widetilde {b}(2k)\left( \sum_{j=1}^n \norm{ a_j }^2 \right)^{\frac{1}{2}}\right)^{2k}
			=
	 \widetilde {b}(2k)^{2k} \left( \sum_{j=1}^n \norm{ a_j }^2 \right)^{k}.
\]
This holds for all $k \in \mathbb{N}$ and therefore for every summand in (\ref{eq:seriesexp}). Note that $\widetilde b (2k)$ does not depend on $n$.

Therefore we get
\begin{align}\label{eq:series_taylor}
		\erw{ \exp \left( \frac{\norm{\sum_{j=1}^n a_j X_j}^{2} }{C^2 \sum_{j=1}^n\norm{ a_j } ^2}\right)}
		&\leq
			\sum_{k=0}^{\infty} \frac{1}{k!} \frac{1}{C^{2k}}\left(\frac 2 N\right) ^k\left( \frac{\Gamma(k+\frac N 2)}{\Gamma (\frac N 2)}\right). 			
\end{align}
Note that $\Gamma(k+\frac N 2)=\Gamma\left( \frac N 2 \right) \prod_{l=1}^k \left(k-l+\frac N 2\right)$.

Consider the function $f(x):=(1-\frac 2 N x)^{-\frac N 2}$. The right-hand side of inequality (\ref{eq:series_taylor}) is the Taylor expansion of the function $f$ at the point $x=\frac 1 {C^2}$.

So we get
\begin{align*}
\erw{ \exp \left( \frac{\norm{\sum_{j=1}^n a_j X_j}^{2} }{C^2 \sum_{j=1}^n \norm{ a_j } ^2}\right)}
&\leq f\left(\frac 1 {C^2}\right).
\end{align*}

Now let $C:=b(N)=\sqrt{\frac 2 N}{\sqrt{\frac{1}{1-(\frac 1 2 )^{\frac 2 N }}}}$. Then $f\left(\frac 1 {C^2}\right)=2$ and this proves that the inequality from Theorem \ref{thm1} holds true.

\section{Proof of the optimality}\label{sec:optimality}
In this section let $X_j$, $j\in \mathbb{N}$ be an i.i.d. family of equidistributed random variables on the sphere $S^{N-1}$. Denote $Y_n:=\sum_{j=1}^n \frac{1}{\sqrt{n}}X_j$.

\begin{lem}\label{lem:uniformint}
Let  $C\geq \sqrt{\frac 2 N}{\sqrt{\frac{1}{1-(\frac 1 2 )^{\frac 2 N }}}}$. Then
the family of random variables
\[
\left(\exp \left (\frac {\norm{\sum_{j=1}^{n}\frac{1}{\sqrt{n}}X_j}}{C} \right)^2 \right), n\in \mathbb{N}
\]
is uniformly integrable.
\begin{proof}
According to \cite[Theorem 6.19]{Klenke2008a} it suffices to prove that for some $p>1$, 
\[
I(p):=\sup _{n\in\mathbb{N}} \erw{\left(\exp \left( \frac {\norm{Y_n}}{C}\right)^2 \right)^p} < \infty.
\]
First note that for a $N$-dimensional Gaussian variable $Z$ we have $\erw {\norm{X_j}^{2k}}=1\leq \erw{\norm{Z}^{2k}}$. Using a theorem of Zolotarev \cite[Theorem 3]{Zolotarev1968a} this implies
\[
\prob{\norm{Y_n}>t}\leq \exp (-Nq(t)),
\]
where $q(t)=\frac 1 2 (t^2-\ln t - 1)$. For large $t$ we have $t^2-\ln t - 1> \gamma t^2$ for some $\gamma$ close to $1$, say $\gamma \in (\frac 1 2,1)$.

Therefore we find
\begin{align*}
I(p)&= \sup_{n\in \mathbb{N}} \int_0^{\infty} \prob{\exp \left(p\frac {\norm{Y_n}^2}{C^2}\right)>t}\D t\\
&=
1+\sup_{n\in \mathbb{N}} \int_1^{\infty} \prob{\norm{Y_n}>\frac{C}{\sqrt{p}}\sqrt{\ln(t)}}\D t\\
&\leq 1+  \int_1^{\infty} t^{-\frac{N}{2}\frac{C^2\gamma}{p}}\D t.
\end{align*}
So we can choose $p\in (1, \frac{N}{2}{C^2\gamma})$ such that the latter integral is finite.
\end{proof}
\end{lem}

\begin{lem}
Let $Z$ be a $N$-dimensional Gaussian variable. Then we have
\[
\norm{Z}_{\psi_{2}}=\frac{\sqrt{2}}{\sqrt{1-(\frac 1 2)^{\frac 2 N}}}.
\] 
\begin{proof}
Let $C>\sqrt{2}$.
We compute
\begin{align*}
\erw{\exp \left( \frac {\norm{Z}^2}{C^2}\right) }
&= 
\frac{1}{(2\pi)^{N/2}} \int_{\mathbb{R}^N} \exp \left( \frac {\norm{x}^2}{C^2}\right) \exp \left(-\frac {\norm{x}^2}{2}\right) \D x\\
&=
\frac{1}{(2\pi)^{N/2}} \int_{\mathbb{R}^N} \exp \left( -\sum_{j=1}^{N} x_j ^2\left(\frac 1 2 -\frac 1 {C^2}\right)\right) \D x \\
&=
\prod_{j=1}^{N} \frac 1{\sqrt{2\pi}}\int_{\mathbb{R}}\exp\left(-\frac 1 2 t^2 \left(\frac {C^2-2}{C^2}\right)\right)\D t \\
&=
\left(\frac {C^2}{C^2-2} \right)^{\frac N 2}.
\end{align*}
Now we have $\left(\frac {C^2}{C^2-2} \right)^{\frac N 2}\leq 2$ if and only if 
$C\geq \sqrt{\frac{2}{1-(\frac 1 2 )^{\frac 2 N}}}$.
This proves the lemma.
\end{proof}
\end{lem}

\begin{lem}
Let $Z$ be a $N$-dimensional Gaussian variable. Then we have
\[
\lim_{n\to\infty}\norm{\sum_{j=1}^n \frac 1 {\sqrt{n}} X_j}_{\psi_2}=\norm{Z}_{\psi_2}.
\]

\begin{proof}
Assume $\limsup_{n\to\infty} \norm{Y_n}_{\psi_2} > \norm{Z}_{\psi_2}$. Then there exists a subsequence $n_k,k\in \mathbb{N}$ and some $\epsilon>0$ such that
\[
\norm{Y_{n_k}}_{\psi_2} > \norm{Z}_{\psi_2} + \epsilon.
\]
According to Lemma \ref{lem:uniformint} the family $\Bigg(\exp \bigg (\frac {\norm{Y_n}}{\norm{Z}_{\psi_2} + \epsilon} \bigg)^2 \Bigg), n\in \mathbb{N}$ is uniformly integrable. Also 
\[
G_n:=\exp \bigg (\frac {\norm{Y_n}}{\norm{Z}_{\psi_2} + \epsilon} \bigg)^2-\exp \bigg (\frac {\norm{Z}}{\norm{Z}_{\psi_2} + \epsilon} \bigg)^2, n\in\mathbb{N}
\]
is uniformly integrable. 
For $M>0$ we have
\[
\int G_n \dprob \leq  \int_{\{G_n \leq M\}}G_n \dprob + \sup_{n\in\mathbb{N}}\int_{\{G_n > M\}}G_n \dprob.
\]
For every fixed $M>0$,  the first integral tends to $0$ for $n\to \infty$ by the central limit theorem.
The second integral tends to $0$ for $M\to \infty$ due to the uniform integrability.
Therefore 
\[
\lim_{n\to\infty}\int \exp\left(\frac {\norm{Y_n}}{\norm{Z}_{\psi_2}+\epsilon}\right)^2\dprob
=\int \exp\left(\frac {\norm{Z}}{\norm{Z}_{\psi_2}+\epsilon}\right)^2\dprob.
\]
This implies
\begin{align*}
2 &\geq \int \exp\left(\frac {\norm{Z}}{\norm{Z}_{\psi_2}}\right) \dprob \\
&>\int \exp\left(\frac {\norm{Z}}{\norm{Z}_{\psi_2}+\epsilon}\right) \dprob \\
&=\lim_{n\to\infty}\int \exp\left(\frac {\norm{Y_n}}{\norm{Z}_{\psi_2}+\epsilon}\right) \dprob \\
&=\lim_{k\to\infty}\int \exp\left(\frac {\norm{Y_{n_k}}}{\norm{Z}_{\psi_2}+\epsilon}\right) \dprob \\
&\geq 2,
\end{align*}
which is a contradiction.
Therefore
$\limsup_{n\to\infty} \norm{Y_n}_{\psi_2}\leq \norm{Z}_{\psi_2}$.
In the same way we show $\liminf_{n\to\infty} \norm{Y_n}_{\psi_2}\geq \norm{Z}_{\psi_2}$.
\end{proof}
\end{lem}
This finishes the proof our Theorem.
\subsection*{Acknowledgment}
I thank Prof. H. K\"onig for his support and advice during my PhD studies.
Part of my research was funded by DFG project KO 962/10-1.

\providecommand{\bysame}{\leavevmode\hbox to3em{\hrulefill}\thinspace}
\providecommand{\MR}{\relax\ifhmode\unskip\space\fi MR }
\providecommand{\MRhref}[2]{%
  \href{http://www.ams.org/mathscinet-getitem?mr=#1}{#2}
}
\providecommand{\href}[2]{#2}

\end{document}